\documentclass[10pt]{article}
\usepackage{amsmath}
\usepackage{amsfonts}
\usepackage{amssymb}
\usepackage{comment}
\usepackage{graphicx,color}
\long\def\proof#1{\removelastskip\vskip\baselineskip\relax\noindent{\it
Proof\if!#1!\else\ \ignorespaces#1\fi.\ }\ignorespaces}

\newcommand{\lgs}[2]{\mbox{$\left(\frac{#1}{#2}\right)$}}
\newcommand{\leg}[2]{\mbox{$\left(\dfrac{#1}{#2}\right)$}}

\newcommand{\ov}{\overline}
\newcommand{\Q}{{\mathbb Q}}
\newcommand{\Z}{{\mathbb Z}}
\newcommand{\R}{{\mathbb R}}

\newcommand{\C}{{\mathbb C}}

\newcommand{\z}{\zeta}
\renewcommand{\th}{\theta}

\newcommand{\G}{\Gamma}

\newcommand{\f}{{\mathfrak f}}
\newcommand{\al}{\alpha}

\newcommand{\ga}{\gamma}

\newcommand{\gd}{{\mathfrak d}}

\newcommand{\eps}{\varepsilon}

\renewcommand{\H}{\mathcal H}
\newcommand{\bftau}{\pmb{\tau}}
\newcommand{\bfz}{\mathbf z}

\DeclareMathOperator{\Tr}{Tr}
\DeclareMathOperator{\NO}{{\mathcal N}}
\DeclareMathOperator{\N}{{\NO\!}}
\DeclareMathOperator{\sign}{sign}

\DeclareMathOperator{\GL}{GL}

\DeclareMathOperator{\SL}{SL}

\DeclareMathOperator{\lcm}{lcm}

\newcommand{\psmm}[4]{\left(\begin{smallmatrix}{#1}&{#2}\\{#3}&{#4}\end{smallmatrix}\right)}

\newcommand{\Proof}{{\it Proof. \/}}
\newcommand{\squareforqed}{\hbox{\rlap{$\sqcap$}$\sqcup$}}
\newcommand{\qed}{\ifmmode\squareforqed\else{\unskip\nobreak\hfil
\penalty50\hskip1em\null\nobreak\hfil\squareforqed
\parfillskip=0pt\finalhyphendemerits=0\endgraf}\fi}

\newcommand{\fp}{\qed\removelastskip\vskip\baselineskip\relax}

\newtheorem{theorem}{Theorem}[section]
\newtheorem{corollary}[theorem]{Corollary}
\newtheorem{proposition}[theorem]{Proposition}
\newtheorem{lemma}[theorem]{Lemma}
\newtheorem{definition}[theorem]{Definition}
\newtheorem{conjecture}[theorem]{Conjecture}

\newcommand{\litem}{\par\noindent\dimen0=\parindent%
    \advance\dimen0 by-4pt
               \hangindent=\dimen0\ltextindent}

\newcommand{\ltextindent}[1]{\hbox to \hangindent{#1\hss}\ignorespaces}
\newcommand{\ltextjndent}[1]{\hbox to \hangindent{#1\hss}\ignorespaces\kern-1ex}

\renewcommand{\pmod}[1]{\allowbreak\ ({\rm{mod}}\,\,#1)}

\newtheorem{remarks}[theorem]{Remarks}

\begin{document}
\pagestyle{plain}
\title{The Shimura Lift and Hilbert Modular Forms}

\author{Henri Cohen}

\maketitle

\begin{abstract}
  After recalling basic formulas concerning modular forms of half-integral
  weight, Gegenbauer polynomials, and Rankin--Cohen brackets, we give four
  versions of the Shimura lift and several applications to Hilbert modular
  forms associated to real quadratic fields. We claim no originality, but
  wanted to collect all the formulas in a single paper.
\end{abstract}  

\section{Introduction}

It is well-known and easy to see that there are close relations between the
Shimura correspondence sending modular forms of half-integral weight
$k+1/2$ to modular forms of integral weight $2k$ and Hilbert modular
forms on real quadratic fields. The aim of the present paper, which does
not claim any originality, is to write down explicitly the precise formulas
and applications. Some parts of this paper are taken from the author's
50-year-old thesis whose title ``modular forms in one and two variables''
is very similar. We refer to the old papers of the author
\cite{Coh1}, \cite{Coh2}, and \cite{Coh3}, to the more recent paper
\cite{Coh4} for many more details on the subjects treated here, as well as
to the two books on modular forms coauthored by the author \cite{Coh-Str}
and \cite{Coh-Gan}.

\smallskip

Since most results are classical, I would be grateful for
corrections and/or improvements to the statements given here.

\smallskip

It is important from the start to fix notation. Whenever possible, we will
use the letters $f$ and $g$, possibly subscripted, to denote
one-variable modular forms of \emph{integral} weight. We will use capital
letters such as $F$ and $G$ to denote modular forms of
\emph{several variables}, with the exception of Eisenstein series which
we will denote by $E_k$ or $F_k$ (see definitions below), but the context
will make the meaning clear. Finally, we will use the letter $h$ possibly
subscripted to denote one-variable modular forms of \emph{half-integral}
weight.

\smallskip

To make everything explicit and essentially self-contained,
we begin by recalling notation for characters, and basic results on modular
forms of half-integral weight, Gegenbauer polynomials, and
Rankin--Cohen brackets.

\section{Preliminary Notation and Results}

\smallskip

\subsection{Notation for Characters}

\smallskip

We will reserve the Greek letter $\chi$ (always subscripted) to denote
quadratic characters: $\chi_D$ is the character defined by
$\chi_D(n)=\lgs{D}{n}$, always with $D\equiv0,1\pmod4$
but not necessarily fundamental. Thus, arbitrary characters will be denoted
by $\psi$, possibly subscripted, and when we say that $\psi$ is a character
modulo $N$ this implies that $\psi(d)=0$ if and only if $\gcd(d,N)>1$.

\smallskip

If $\psi$ is a general Dirichlet character modulo $N$ we denote by
$\f\mid N$ the conductor of $\psi$ (not to be confused with a function $f$),
and by $\psi_{\f}$ the character modulo $\f$ equivalent to $\psi$. More
generally, for any $M$ such that $\f\mid M\mid N$ we will denote by $\psi_M$
the character modulo $M$ equivalent to $\psi$, defined by $\psi_M(m)=0$ if
$\gcd(m,M)>1$ and $\psi_M(m)=\psi_{\f}(m)$ if $\gcd(m,M)=1$. Otherwise said,
if we denote by $\chi_{0,M}$ the trivial character modulo $M$ (the indicator
function of integers coprime to $M$), we have $\psi_M=\chi_{0,M}\psi_f$.
(note that if we wanted to be completely consistent, using our notation
for quadratic characters we should write $\chi_{M^2}$ instead of $\chi_{0,M}$).

\smallskip

Finally, we will very often use the \emph{product} of a character $\psi$ and
a character $\chi_D$: if we write $\psi\chi_D$ this means that
$\psi\chi_D(m)=\psi(m)\chi_D(m)$, in other words
we do \emph{not} take the primitive character equivalent to it.
If $\psi$ is a character modulo $N$ the character $\psi\chi_D$ is a character
modulo $\lcm(N,|D|)$. For instance $\chi_D\chi_D$ is the trivial
character $\chi_{0,|D|}$ modulo $|D|$, and \emph{not} the trivial
character $\chi_1$ (and indeed, $\chi_D\chi_D=\chi_{D^2}$; see the above
note).

\smallskip

\subsection{Modular Forms of Half-Integral Weight}

\smallskip

Recall that we denote by
$\th(\tau)=\sum_{s\in\Z}q^{s^2}=1+2\sum_{n\ge1}q^{n^2}$ the standard
one-variable theta function of weight $1/2$ on $\G_0(4)$ satisfying
$\th(\ga(\tau))=v_{\th}(\ga)(c\tau+d)^{1/2}\th(\tau)$ for
$\ga=\psmm{a}{b}{c}{d}\in\G_0(4)$, with $v_{\th}(\ga)=\lgs{-4}{d}^{-1/2}\lgs{c}{d}$.

\smallskip

Let $h=\sum_{n\ge0}a(n)q^n\in M_{k+1/2}(\G_0(N),\psi)$ be a modular form
of half-integral weight. This notation, which I believe is due to Shimura, is not
really appropriate because of the multiplier system.
In addition to the holomorphy conditions on $\H$ and at the
cusps, this means that $g=h/\th^{2k+1}$ is a modular function (of weight $0$)
and character $\psi$ on $\G_0(N)$, in other words that $g|_0\ga=\psi(d)g$
for $\ga=\psmm{a}{b}{c}{d}\in\G_0(N)$. If $M$ is the level of the modular
function $g$ the level of $h$ is $\lcm(M,4)$, hence the level of
$h=g\th^{2k+1}$ is also $M=\lcm(M,4)$, in other words $M=N$ has to be
divisible by $4$. The modularity condition is that
$h|_{k+1/2}\ga=\psi(d)v_{\th}(\ga)^{2k+1}h$.
Now $v_{\th}(\ga)^{2k}=\lgs{-4}{d}^k$, so $h|_{k+1/2}\ga=(\chi_{-4}^k\psi)(d)v_{\th}(\ga)h$. When $k$ is even, the notation $h\in M_{k+1/2}(\G_0(N),\psi)$ is
justified, but when $k$ is odd, one should write instead
$h\in M_{k+1/2}(\G_0(N),\chi_{-4}\psi)$, in both cases keeping the multiplier
system $v_{\th}$ (instead of $v_{\th}^{2k+1}$) implicit. Looking at what
happens when changing $\ga$ to $-\ga$ shows that to have a nonzero theory,
when $k$ is even $\psi$ must be an even character, and when $k$ is odd the
character must be an odd character, but since this is $\chi_{-4}\psi$, this
again means that $\psi$ must be an even character.
Thus, since we cannot change established notation, when
writing $h\in M_{k+1/2}(\G_0(N),\psi)$ we always assume that $4\mid N$
and $\psi$ is an even character. Nonetheless, even if trivially the same,
when $k$ is odd it would be better to say that $\chi_{-4}\psi$ is an odd
character than to say that $\psi$ is an even one.

\smallskip

Another consequence of the notation is the following not completely intuitive
but easy to prove result, which should be considered as a warning:

\begin{proposition}\label{prop:half}
 Let $f\in M_k(\G_0(N),\psi)$, and $M$ a positive integer.
\begin{enumerate}\item We have $f(M\tau)\in M_k(\G_0(NM),\psi\psi_1)$,
where $\psi_1=1$ if $k\in\Z$ and $\psi_1=\chi_{4M}$ if $k\in 1/2+\Z$.
\item Let $f_i\in M_{k_i}(\G_0(N_i),\psi_i)$. We have
$f_1f_2\in M_{k_1+k_2}(\G_0(\lcm(N_1,N_2)),\psi_1\psi_2\psi_{1,2})$, with
$\psi_{1,2}=1$ if $k_1$ and $k_2$ are in $\Z$, $\psi_{1,2}=\chi_{-4}^{k_i}$
if $k_i\in\Z$ and $k_{3-i}\in1/2+\Z$, and $\psi_{1,2}=\chi_{-4}^{k_1+k_2}$ if
$k_1$ and $k_2$ are in $1/2+\Z$.
\item In particular, we have $f\th\in M_{k+1/2}(\G_0(\lcm(N,4)),\psi\psi_1)$,
where $\psi_1=\chi_{-4}^k$ if $k\in\Z$, and $\psi_1=\chi_{-4}^{k+1/2}$
if $k\in1/2+\Z$.
\end{enumerate}\end{proposition}

Note in particular in (1) the appearance of a quadratic character other than
$\chi_{-4}$.

Results (2) and (3) of the above proposition are more generally valid
if we replace the products by the $\ell$-th Rankin--Cohen brackets (see below),
of course after increasing the weight by $2\ell$.

\smallskip

\subsection{Variations on Gegenbauer Polynomials}\label{sec:geg}

\smallskip

Recall the following:

\begin{definition} We define the Gegenbauer polynomials by the generating
  function
  $$\dfrac{1}{(1-2xt+t^2)^a}=\sum_{\ell\ge0}C_{\ell}^{(a)}(x)t^\ell$$
\end{definition}

In number theory, these polynomials occur naturally for instance in the
Eichler--Selberg trace formula. In applications of Rankin--Cohen brackets
we will use a different family of polynomials which are closely related:

\begin{definition}
For $\ell$ a nonnegative integer, we define the polynomials
$$P_{k,\ell}(n,s)=\sum_{j=0}^{\ell}(-1)^j\binom{\ell+k-1}{j}\binom{\ell-1/2}{\ell-j}(n-s^2)^{\ell-j}s^{2j}\;.$$
\end{definition}

For instance, we have
\begin{align*}
P_{k,0}(n,s)&=1\;,\\
P_{k,1}(n,s)&=n/2-(1/2)(2k+1)s^2\;,\\
P_{k,2}(n,s)&=3n^2/8-(3/4)(2k+3)ns^2+(1/8)(2k+3)(2k+5)s^4\;.
\end{align*}

The following result shows the relation with Gegenbauer polynomials and
additional identities. The proofs are combinatorial and not difficult,
and are left to the reader.

\begin{proposition}\label{prop:pkl}
  \begin{enumerate}
  \item We have
    $$P_{k,\ell}(n,s)=(-1)^\ell\frac{\displaystyle\binom{2\ell}{\ell}}{\displaystyle\binom{k+\ell-3/2}{\ell}}(n/4)^{\ell}C_{2\ell}^{(k-1/2)}\left(\dfrac{s}{\sqrt{n}}\right)\;.$$
  \item We have
    $$P_{k,\ell}(n,s)=\binom{2\ell}{\ell}\sum_{j=0}^{\ell}(-1)^j\binom{k+\ell-3/2+j}{j}\frac{\displaystyle\binom{\ell}{j}}{\displaystyle\binom{2j}{j}}(n/4)^{\ell-j}s^{2j}\;.$$
  \item We have
    $$P_{k,\ell}(n^2D,s)=(-1)^\ell\dfrac{(2\ell)!}{2^{2\ell}\ell!}\dfrac{(k+2\ell-1)!(k+\ell-1)!}{(k-1)!^2}D^\ell Q_{k,2\ell}(n,s,D)\;,$$
    where the polynomials $Q_{k,\ell}$ will be defined in Corollary
    \ref{cor:qk} below.
  \end{enumerate}
\end{proposition}

\smallskip

\subsection{Rankin--Cohen Brackets}

\smallskip

Recall that for any $\ga=\psmm{a}{b}{c}{d}\in\GL_2^+(\R)$ and
any function $f$ we define $f|_k\ga(\tau)=(ad-bc)^{k/2}(c\tau+d)^{-k}f((a\tau+b)/(c\tau+d))$. For notational simplicity, we write $E=E(\tau)=(\tau-\ov{\tau})^{-1}=1/(2i\Im(\tau))$, and we will denote by $\partial_{\tau}$ the derivative with
respect to $\tau$. We need some preliminary results.

\begin{proposition} Let $f$ be any $C^\infty$ function and
$\ga=\psmm{a}{b}{c}{d}\in\GL_2^+(\R)$.
\begin{enumerate}
\item We have $E|_2\ga=E-c/(c\tau+d)$ and $\partial_\tau(E)=-E^2$.
\item We have $(\partial_\tau f+kEf)|_{k+2}(\ga)=\partial_\tau(f|_k\ga)+kE(f|_k\ga)$.
\item For a nonnegative integer $r$ define
$$A_r(f)=\sum_{0\le m\le r}\binom{r}{m}\binom{k-1+r}{r-m}\dfrac{E^{r-m}}{(r-m)!}\partial_\tau^m(f)\;.$$
We have $A_{r+1}(f)=\partial_\tau(A_r(f))+(k+2r)EA_r(f)$ and
$A_r(f)|_{k+2r}\ga=A_r(f|_k\ga)$.
\end{enumerate}\end{proposition}

\Proof (1) and (2), and the first part of (3), are simple verifications,
and the second part of (3) follows from the first by induction on $r$.\fp

\begin{corollary}\label{cor:FG}
Define the formal power series in $z$
$$A_f(\tau,z)=\sum_{r\ge0}\dfrac{(2\pi iz^2)^rA_r(f)(\tau)}{r!(k-1+r)!}\text{\quad and\quad}B_f(\tau,z)=\sum_{r\ge0}\dfrac{(2\pi iz^2)^r\partial_\tau^r(f)(\tau)}{r!(k-1+r)!}\;,$$
and let $\ga=\psmm{a}{b}{c}{d}\in GL_2^+(\R)$ act by
$$A_f(\tau,z)|_k\ga=(ad-bc)^{k/2}(c\tau+d)^{-k}A_f\left(\dfrac{a\tau+b}{c\tau+d},\dfrac{(ad-bc)^{1/2}}{c\tau+d}z\right)\;,$$
and similarly for $B_f$. We have
$A_f(\tau,z)=e^{2\pi iz^2E(\tau)}B_f(\tau,z)$ and
$$A_f|_k\ga(\tau,z)=A_{f|_k\ga}(\tau,z)\text{\quad and\quad}
B_f|_k\ga(\tau,z)=e^{2\pi iz^2c/(c\tau+d)}B_{f|_k\ga}(\tau,z)\;.$$
\end{corollary}

\Proof The first statement is immediate by series manipulation,
and the other statements follow immediately from the proposition and the relation
between $A_f$ and $B_f$.\fp

\begin{proposition} For $i=1,2$ let $f_i\in M_{k_i}(\G_0(N_i),\psi_i)$ with
  $k_i$ integral or half-integral, and for $\ell\ge0$ define
  $$[f_1,f_2]_{\ell}=\dfrac{1}{(2\pi i)^\ell}\sum_{0\le j\le \ell}(-1)^j\dfrac{(k_1+\ell-1)!(k_2+\ell-1)!}{(k_1+j-1)!(k_2+\ell-j-1)!}\dfrac{\partial_{\tau}^j(f_1)\partial_{\tau}^{\ell-j}(f_2)}{j!(\ell-j)!}\;.$$
  We have $[f_1,f_2]_{\ell}\in M_{k_1+k_2+2\ell}(\G_0(\lcm(N_1,N_2)),\psi_1\psi_2\psi_{1,2})$,
where $\psi_{1,2}$ is given in Proposition \ref{prop:half} (2),
  and it is a cusp form if $\ell>0$.
\end{proposition}

\Proof Set $C(\tau,z)=B_{f_1}(\tau,z)B_{f_2}(\tau,iz)$. By the above
corollary we know that $C|_{k_1+k_2}\ga(\tau,z)=\psi_1\psi_2\psi_{1,2}(d)C(\tau,z)$,
and since $C(\tau,z)$ is a formal power series in $z$ (in fact in $z^2$),
and $z$ is acted upon by $z\mapsto z/(c\tau+d)$, it follows that the
coefficient of $z^{2\ell}$ belongs to
$M_{k_1+k_2+2\ell}(\G_0(N),\psi_1\psi_2\psi_{1,2})$, and a small
computation shows that up to a multiplicative constant it is equal to
$[f_1,f_2]_{\ell}$.

For the last statement, note that the above argument shows more generally
that for all $\ga\in\SL_2(\Z)$ we have
$([f_1,f_2]_{\ell})|_{k_1+k_2+2\ell}\ga=[f_1|_{k_1}\ga,f_2|_{k_2}\ga]_{\ell}$,
hence when $\ell>0$ since at least one derivative kills the constant term,
the bracket vanishes at all cusps.\fp

\smallskip

We will also need the following generalization of the above proposition
where one of the $f_i$ is equal to the quasi-modular function $E_2$.
More precisely, we have the following result whose easy proof is left to
the reader:
\begin{lemma}\label{lem:E2}
  Let $f\in M_k(\G_0(N),\psi)$ with $k$ integral or half-integral, and define
  $$[f,E_2(M\tau)]^*_\ell=[f,E_2(M\tau)]_\ell-\dfrac{12\ell!}{M(k+\ell)}\dfrac{\partial_{\tau}^{\ell+1}(f)}{(2\pi i)^{\ell+1}}\;.$$
  We have $[f,E_2(M\tau)]^*_\ell\in M_{k+2+2\ell}(\G_0(\lcm(M,N)),\psi)$,
  and it is a cusp form if $\ell>0$.
\end{lemma}

\smallskip

\subsection{An Application of RC brackets}

\smallskip

\begin{proposition}\label{prop:hlm} Let $M$ be a positive integer and $\ell$
 a nonnegative integer.
  \begin{enumerate}
  \item Let $f=\sum_{n\ge0}a(n)q^n\in M_k(\G_0(N),\psi)$ with $k$ integral
or half-integral, and set
  $$A_{\ell,M}(n)=\sum_{\substack{|s|\le n^{1/2}\\s^2\equiv n\pmod M}}P_{k,\ell}(n,s)a\left(\dfrac{n-s^2}{M}\right)\text{\ and\ }h_{\ell,M}=\sum_{n\ge0}A_{\ell,M}(n)q^n\;,$$
where the polynomials $P_{k,\ell}$ have been defined in Section \ref{sec:geg}.

We have $h_{\ell,M}\in M_{k+2\ell+1/2}(\G_0(\lcm(4,MN)),\psi\psi_1)$,
where $\psi_1=\chi_{-4}^k$ if $k\in\Z$, $\psi_1=\chi_{-4}^{k+1/2}$ if
$k\in1/2+\Z$, and it is a cusp form if $\ell>0$.
\item The same result is true for $f=-E_2/24$, so that $a(n)=\sigma_1(n)$
  for $n\ge1$ and $a(0)=-1/24$, with $k=2$, $N=1$, and $\psi=1$, if we define
    $$h_{\ell,M}=\sum_{n\ge0}A_{\ell,M}(n)q^n+\dfrac{(-1)^\ell}{(2\ell+1)M}\sum_{s\in\Z}s^{2\ell+2}q^{s^2}\;.$$
    \end{enumerate}
\end{proposition}
Note that $A_{\ell,M}(0)=\delta_{\ell,0}a(0)$.

\smallskip

\Proof Indeed, for (1) it is immediate to check that
$h_{\ell,M}(\tau)$ is a constant multiple of $[\th(\tau),f(M\tau)]_\ell$, and
for (2) that it is a constant multiple of 
$-(1/24)[\th(\tau),E_2(M\tau)]^*_{\ell}$.\fp

\section{Variations on the Shimura Correspondence}

There are (at least) four versions of the Shimura lift from modular
forms of half integral weight $k+1/2$ to modular forms of integral weight
$2k$, differing only at the prime $p=2$. The functorial way to do this,
initiated by Shimura, is to define it for eigenforms by giving the $L$-function
of the lift. We will instead use a less elegant way initiated probably by
Zagier and the author, by working directly on $q$-expansions of functions
which are not necessarily eigenforms. To make things completely transparent,
we will give more details than usual.

\subsection{The Four Variations}

\medskip

In this section we fix $h=\sum_{n\ge0}c(n)q^n\in M_{k+1/2}(\G_0(N),\psi)$
with $k\ge1$ integral. We use on purpose the notation $c(n)$ instead of $a(n)$
to insist that it is the coefficient of a modular form of half-integral
weight.

\smallskip

The first definition of the lift is essentially Shimura's initial definition,
completed by works of Cipra and Niwa, see \cite{Shi}, \cite{Niw}, \cite{Cip}:

\begin{definition}\label{def:shi0} Let $t$ be a positive squarefree
  integer, and let $D=(-1)^kt$ or $D=(-1)^k4t$ be the fundamental discriminant
  corresponding to $(-1)^kt$. For $n\ge1$, set
  $$C_0(n;t)=\sum_{d\mid n}\psi\chi_D(d)d^{k-1}c((n/d)^2t)\;.$$
  We define the Shimura lift of $h$ at $t$ by
  $$S_0(h;t)=\dfrac{c(0)}{2}L(\psi\chi_D,1-k)+\sum_{n\ge1}C_0(n;t)q^n\;.$$
\end{definition}

\smallskip

\begin{remarks}\begin{enumerate}
  \item When we speak of fundamental discriminants $D$, we always include
    $D=1$.
  \item Even if its conductor is smaller than $N$, it is important to
    consider $\psi$ as a character modulo $N$, and hence $\psi\chi_D$
    as a character modulo $\lcm(|D|,N)$, i.e., vanishing on integers not
    coprime to $\lcm(D,N)$. For instance, if $\psi$ is trivial
    (no nebentypus), it must be considered as the trivial character
    $\chi_{0,N}$ modulo $N$.
  \item In particular, since $L(\psi\chi_D,1-k)$ is the value at $1-k$ of a
    not necessarily primitive character, beware that
    some software systems such as {\tt Pari/GP} always return the value of the
    corresponding primitive character, so it is essential to multiply
    by suitable Euler factors.
\end{enumerate}\end{remarks}

\smallskip

In this context, a small part of Shimura's theorem completed by Niwa and Cipra
is as follows:
\begin{theorem}\label{thmshi0}
  If $k\ge1$ we have $S_0(h;t)\in M_{2k}(\G_0(N/2),\psi^2)$.
\end{theorem}

Note that there are numerous additional results giving conditions under which
$S_0(h;t)$ is a cusp form when $h$ is, compatibility with Hecke operators,
etc., which do not concern us in the present paper.

\smallskip

The second definition is a variant of the first and lies between the first and
the third definition, which is essential for using Kohnen's theorem.

\begin{definition}\label{def:shi1} Let $D$ be a fundamental discriminant such
  that $\sign(D)=(-1)^k$. For $n\ge1$, set
  $$C_1(n;D)=\sum_{d\mid n}\psi\chi_D(d)d^{k-1}c((n/d)^2|D|)\;.$$
  We define the Shimura lift of $h$ at $D$ by
  $$S_1(h;D)=\dfrac{c(0)}{2}L(\psi\chi_D,1-k)+\sum_{n\ge1}C_1(n;D)q^n\;.$$
\end{definition}

When $D\equiv1\pmod4$ this is identical to the previous definition, but
when $D\equiv0\pmod4$ it is slightly different (see below for the comparison).
We also have:
\begin{theorem}\label{thmshi1}
  If $k\ge1$ we have $S_1(h;D)\in M_{2k}(\G_0(N/2),\psi^2)$.
\end{theorem}

\smallskip

The third definition is more or less the same as the previous one,
but introduces important restrictions needed to state Kohnen's theorem; see \cite{Koh1}, \cite{Koh2}, and
\cite{Koh-Zag}:

\begin{definition}\label{def:shi2} Let $D$ be a fundamental discriminant such
  that $\sign(D)=(-1)^k$. Assume that the conductor $\f$ of $\psi$ divides
  $N/4$, and as explained above, denote by $\psi_{N/4}$ the character modulo
  $N/4$ equivalent to $\psi$. For $n\ge1$, set
  $$C_2(n;D)=\sum_{d\mid n}\psi_{N/4}\chi_D(d)d^{k-1}c((n/d)^2|D|)\;.$$
  We define the Shimura--Kohnen lift of $h$ at $D$ by
  $$S_2(h;D)=\dfrac{c(0)}{2}L(\psi_{N/4}\chi_D,1-k)+\sum_{n\ge1}C_2(n;D)q^n\;.$$
\end{definition}

The same remark applies as before to the character being defined modulo $\lcm(N/4,|D|)$.
In particular, the divisor sum is over divisors coprime to $N/4$ (and to $D$
of course), while in the previous definition the sum was over divisors
coprime to $N$.

\begin{theorem}\label{thmshi2} Assume as always that $k\ge1$.
  \begin{enumerate}
  \item We have $S_2(h;D)\in M_{2k}(\G_0(N/2),\psi^2)$.
  \item If $h$ belongs to the Kohnen $+$-space, in other words if
    $c(n)=0$ when $(-1)^kn\equiv2,3\pmod4$, then if in addition
    $N/4$ is squarefree, the conductor of $\psi$ divides $N/4$, and $\psi$
    is a quadratic character, we have
    $S_2(h;D)\in M_{2k}(\G_0(N/4),\psi^2)=M_{2k}(\G_0(N/4))$.
  \end{enumerate}
\end{theorem}
This is proved in loc.~cit. when $N/4$ is odd and squarefree, and the
generalization to all squarefree $N/4$ is in \cite{Man} and \cite{MRV}.

Note that this result seems to be true under weaker conditions on $N$ and
$\psi$. More precisely, on the basis of a small number of experiments,
it is possible that, keeping the assumptions on $\psi$, the conclusion
is valid for any $N$, not only with $N/4$ squarefree. I would be interested
to know if one can exhibit a counterexample (or prove this conjecture).

\smallskip

Finally, there is a variant of Definition \ref{def:shi2} as follows,
see \cite{MRV} and \cite{Ued-Yam}:

\begin{definition}\label{def:shi3} Assume that there exists a Dirichlet
  character $\psi'$ modulo $N/4$ (necessarily an odd character) such that
  $\psi=\chi_{-4}\psi'$, and let $D$ be a fundamental
  discriminant congruent to $0$ modulo $4$ such that $\sign(D)=(-1)^{k+1}$.
  For $n\ge1$, set
  $$C_3(n;D)=\sum_{d\mid n}\psi'\chi_D(d)d^{k-1}c((n/d)^2|D|)\;.$$
  We define the Shimura-Kohnen lift of $h$ at $D$ by
  $$S_3(h;D)=\dfrac{c(0)}{2}L(\psi'\chi_D,1-k)+\sum_{n\ge1}C_3(n;D)q^n\;.$$
\end{definition}

Theorem \ref{thmshi2} again holds with $S_2$ replaced by $S_3$.

\smallskip

\begin{remarks}
\begin{enumerate}\item The condition on $\psi$ is automatically satisfied when
  $N/4$ is odd.
\item The restriction to even fundamental discriminants is essential: when
  $D$ is odd, the lift need not even be a modular form.
\end{enumerate}
\end{remarks}

\medskip

\subsection{Relation Between the Definitions}

\medskip

\begin{proposition}\begin{enumerate}
  \item We have
    $$S_1(h;D)(\tau)=S_2(h;D)(\tau)-\psi_{N/4}\chi_D(2)S_2(h;2\tau)\;.$$
    In particular, if $N/4$ is even or $D$ is even we have $S_1(h;D)=S_2(h;D)$.
  \item If $D$ is odd we have
    $$S_1(h;D)=S_0(h;(-1)^kD)=S_0(h;|D|)\;,$$
  \item If $D$ is even we have
    $$S_1(h;D)(\tau)=S_0(h;|D|/4)(\tau)|_{2k}U_2=\dfrac{S_0(h;|D|/4)(\tau/2)+S_0(h;|D|/4)((\tau+1)/2)}{2}\;,$$
    $$S_3(h;D)(\tau)=S_0(h;|D|/4)(\tau)|_{2k}U_2=\dfrac{S_0(h;|D|/4)(\tau/2)+S_0(h;|D|/4)((\tau+1)/2)}{2}\;.$$
    By abuse of notation we can write this as $S_3(h;D)=S_1(h;-D)$, even
    though $-D$ may not be a fundamental discriminant.
  \end{enumerate}
\end{proposition}

\Proof (1). The difference between the definitions of $S_1$ and $S_2$ is that
in $S_2$ we sum over divisors $d$ of $n$ coprime to $N/4$ instead of $N$. Thus
it is clear that we have
$$C_2(n;D)=C_1(n;D)+\sum_{\substack{d\mid n\\\gcd(d,N/4)=1\\\gcd(d,N)>1}}\psi_{N/4}\chi_D(d)d^{k-1}c((n/d)^2|D|)\;.$$
If $N/4$ is even, $\gcd(d,N/4)=1$ implies that $\gcd(d,N)=1$ so the second
sum is empty, hence $S_2(h;D)=S_1(h;D)$. If $N/4$ is odd, the second sum
is thus over even divisors of $n$ coprime to $N/4$, hence setting $d=2d'$
it is clear that $C_2(n;D)=C_1(n;D)+\psi_{N/4}\chi_D(2)2^{k-1}C_2(n/2;D)$
(with $C_2(n/2;D)=0$ if $n$ is odd), hence
$$S_2(h;D)(\tau)=S_1(h;D)(\tau)+\psi_{N/4}\chi_D(2)2^{k-1}S_2(h;D)(2\tau)\;,$$
proving the claim since $\chi_D(2)=0$ when $D$ is even.

\medskip

(2). If $D$ is odd, to compare $S_0$ and $S_1$ we evidently choose
$t=|D|=(-1)^kD$, and in this case clearly $S_0(h;t)=S_1(h;D)$.

\medskip

(3). We now assume that $D$ even, in other words $D\equiv0\pmod4$. To compare
$S_0$ and $S_1$ we must choose $t=|D|/4=(-1)^kD/4$. Since $\psi$ is a
character modulo $N$ and $4\mid N$, we have $\psi\chi_D=\psi\chi_t$. Thus
$$C_1(n;D)=\sum_{d\mid n}\psi\chi_t(d)d^{k-1}c((n/d)^2(4t))\;.$$
Since the even divisors do not contribute, we have
$$C_0(2n;t)=\sum_{d\mid n}\psi\chi_t(d)d^{k-1}c((2n/d)^2t)=C_1(n;D)$$
It follows that
$$S_1(h;D)(\tau)=S_0(h;t)(\tau)|_{2k}U_2=\dfrac{S_0(h;D)(\tau/2)+S_0(h;D)((\tau+1)/2)}{2}\;.$$
A similar reasoning gives the formula for $S_3$.

Note that, although the formulas for $S_1$ and $S_3$ are identical,
that for $S_1$ applies to discriminants of sign $(-1)^k$, while that
for $S_3$ to discriminants of sign $(-1)^{k+1}$, whence the possible
abuse of notation $S_3(h;D)=S_1(h;-D)$.\fp

\section{Application of the Shimura Correspondence}

\smallskip

\subsection{The Main Theorem}

\smallskip

The main result of this section is as follows:

\begin{theorem}\label{thm4} Let
  $f=\sum_{n\ge0}a(n)q^n\in M_k(\G_0(N),\psi)$ be a modular form of
  integral weight $k$ (so that $\psi$ is a character modulo $N$ such that
  $\psi(-1)=(-1)^k$). For any fundamental discriminant $D$ of sign $(-1)^k$ let
  $\psi_1=\psi\chi_D$ if $k$ is even and $\psi_1=\psi\chi_{-4D}$ if
  $k$ is odd, and let $g_{M,D}(\tau)=\sum_{n\ge0}c_{k,M,D}(n)q^n$ with
  \begin{align*}
    c_{k,M,D}(0)&=\delta_{\ell,0}a(0)\dfrac{L(\psi_1,1-k)}{2}\text{\quad and for $n\ge1$}\\
    c_{k,M,D}(n)&=\sum_{d\mid n}\psi_1(d)d^{k+2\ell-1}A_{\ell,M}((n/d)^2|D|)\;,\end{align*}
  where $\delta_{\ell,0}$ is the Kronecker delta and $A_{\ell,M}$ has been
defined in Proposition \ref{prop:hlm}.
  \begin{enumerate}
  \item We have
    \begin{align*}
      g_{4,D}&\in M_{2k+4\ell}(\G_0(2N/\gcd(4,D,N)),\psi^2)\text{\quad and}\\
      g_{1,D}&\in M_{2k+4\ell}(\G_0(2N)/\gcd(4,N),\psi^2)\;,
    \end{align*}
    and they are cusp forms if either $f$ is a cusp form or if $\ell>0$.
    In particular their level divides $2N$.
  \item Assume that $N\equiv0\pmod4$. We have
    $g_{1,D}\in M_{2k+4\ell}(\G_0(N/2),\psi^2)$, and if in addition
    $D\equiv0\pmod4$ we also have
    $g_{4,D}\in M_{2k+4\ell}(\G_0(N/2),\psi^2)$.
  \item Assume that $N\equiv2\pmod4$. We have
    $g_{1,D}\in M_{2k+4\ell}(\G_0(N),\psi^2)$ and
    $g_{4,\delta D}\in M_{2k+4\ell}(\G_0(N),\psi^2)$, where
    $\delta=1$ if $D\equiv0\pmod4$ and $\delta=4$ if $D\equiv1\pmod4$.
  \item If $k$ is even, $N$ is squarefree, and $\psi$ is a quadratic character
    (including the trivial one), we also have
    $g_{4,D}\in M_{2k+4\ell}(\G_0(N),\psi^2)$.
  \item The above results are also true for $f=E_2(N\tau)$ for any
    $N$, even though $f$ is only quasi modular.
\end{enumerate}
\end{theorem}

Note that when $\delta=4$ we use the notation $g_{4,\delta D}$ even
though $\delta D$ is not fundamental.

\smallskip

\Proof By Proposition \ref{prop:hlm} and using its notation, we know that
$$h_{\ell,M}=\delta_{\ell,0}a(0)+\sum_{n\ge1}A_{\ell,M}(n)q^n\in M_{k+2\ell+1/2}(\G_0(\lcm(4,MN)),\psi\chi_{-4}^k)\;.$$

  \smallskip
  
  (1). We first apply Definition \ref{def:shi1} and Theorem \ref{thmshi1} to
  $h_{\ell,M}$, which does not require any Kohnen-type conditions. The character
  denoted $\psi\chi_D$ in that theorem is thus $\psi\chi_D$ if $k$ is even and
  $\psi\chi_{-4D}$ if $k$ is odd, which is what we denote by $\psi_1$.
  Thus, for $n\ge1$ we have
  $$C_1(n;D)=\sum_{d\mid n}\psi_1(d)d^{k+2\ell-1}\sum_{n\ge1}A_{\ell,M}((n/d)^2|D|)\;,$$
  which is equal to $c_{k,M,D}(n)$, and the constant coefficient is
  $\delta_{\ell,0}a(0)L(\psi_1,1-k)/2$.

  Thus Shimura--Niwa--Cipra's Theorem \ref{thmshi1} tells us that
  the lift $\sum_{n\ge0}C_1(n;D)q^n$ is a modular form
  of weight $2k+4\ell$ and Nebentypus $\psi^2$ (since $\chi_{-4D}^2$ is
  trivial), of level (dividing) $\lcm(4,MN)/2$. For $M=1$, this is equal
  to $2N/\lcm(4,N)$, and for $M=4$ it is equal to $2N$, proving (1) for
  $g_{1,D}$, and also for $g_{4,D}$ when $D\equiv1\mod4$ since in that case
  $\gcd(4,D,N)=1$.

  Thus, assume now that $D\equiv0\pmod4$, so that $(-1)^kD/4=t>0$ and
  squarefree. Here, we apply Definition \ref{def:shi0} and Theorem
  \ref{thmshi0} to the form $h_{\ell,1}$. The character $\psi_1$ is the same,
  and we compute once again that
  $$C_0(n;t)=\sum_{d\mid n}\psi_1(d)d^{k+2\ell-1}\sum_{n\ge1}A_{\ell,1}((n/d)^2t)\;,$$
  with $A_{\ell,1}(m)=\sum_{|s|<m^{1/2}}P_{k,\ell}(m,s)a(m-s^2)$. Since
  $D\equiv0\mod4$ the congruence $(n/d)^2|D|\equiv s^2\pmod4$ in $A_{\ell,4}$
  means that $s$ is even, so writing $s=2s'$ we have
  \begin{align*}A_{\ell,4}((n/d)^2|D|)&=\sum_{|s|<((n/d)^2|D|)^{1/2}}P_{k,\ell}((n/d)^2|D|,s)a\left(\dfrac{(n/d)^2|D|-s^2}{4}\right)\\
    &=\sum_{|s'|<((n/d)^2t)^{1/2}}P_{k,\ell}((n/d)^24t,2s')a((n/d)^2t-{s'}^2)=2^{2\ell}C_0(n;t)\end{align*}
  since $P(4m,2s)=2^{2\ell}P(m,s)$. On the other hand, since
  $h_{\ell,1}$ has level $\lcm(4,N)=4N/\gcd(4,N)$, Theorem \ref{thmshi0} tells
 us that the lift has level $2N/\gcd(4,N)=2N/\gcd(4,D,N)$, proving (1).

  \smallskip

  (2) and the statements of (3) for $D\equiv0\pmod4$ are trivial consequences
  of (1). So assume that $N\equiv2\pmod4$ and $D\equiv1\pmod4$. By (1) the
  form $g_{4,D}$ has level dividing $2N$ which is divisible by $4$, so
  applying the Hecke operator $U_2$ we deduce that $g_{4,D}|U_2$ has level
  dividing $N$, and since $\psi_1(d)=0$ when $2\mid d$ it is clear that
  $g_{4,D}|U_2=g_{4,4D}$ (even though $4D$ is not fundamental), proving (3).
  
  \smallskip

  (4). We now want to apply Kohnen's theorem and generalizations.
  Here (unless $f$ has special properties), we must have $k$ even,
  and then clearly $[\th(\tau),f(4\tau)]_\ell$ belongs to
  the Kohnen $+$-space (unfortunately when $k$ is odd, the vanishing
  coefficients are not the ones needed for Kohnen, and there is apparently
  no way to correct this) in level $\lcm(4,4N)=4N$ with nebentypus
  $\psi$ since $k$ is even. Thus, to apply Definition \ref{def:shi2}
  and Theorem \ref{thmshi2} we must now assume that $N$ is squarefree and
  $\psi$ is a quadratic character, necessarily even, and which automatically
  has conductor dividing $N$ since it comes from a form $f$ of level $N$,
  so equal to what is denoted $\psi_{N/4}$ in Definition \ref{def:shi2}.
  Thus we again have $C_2(n;D)=c_{k,4,D}(n)$ (and the same constant term),
  proving (4).

  \smallskip
  
  (5). By Proposition \ref{prop:hlm} (2), we must add to the definition of
$h_{\ell,M}$ a constant multiple of $\sum_{s\in\Z}s^{2\ell+2}q^{s^2}$.
We apply the same reasoning as above, so if we denote by $b^*$ the Fourier
  coefficients of this modified version and by $b(m)$ those of the initial
 version, we have $b^*(m)=b(m)$ is $m$ is not a square, 
and $b^*(m)=b(m)+Cm^{\ell+1}$ for some constant $C$ if
  $m=s^2$ is a square. Since $k=2$ is even we have $D>0$, hence $|D|$ is not
  a square (this could have happened for $k$ odd for $D=-4$), so that
  $b^*((n/d)^2|D|)=b((n/d)^2|D|)$, proving the result.
\fp

\begin{definition}\label{def:lev} We define a quantity $L(N,D)$ as follows.
  \begin{enumerate}\item If $N\equiv D\equiv0\pmod4$, $L(N,D)=N/2$.
  \item If $N\equiv2\pmod4$ and $D\equiv0\pmod4$, or if $k$ is even,
    $N$ is squarefree, and $\psi$ is a quadratic character, $L(N,D)=N$.
  \item Otherwise $L(N,D)=2N$.
  \end{enumerate}
\end{definition}

Note that $L(N,D)$ depends also on $k$ and $\psi$.

We can then summarize part of the above theorem by saying that
$g_{4,D}\in M_{2k+4\ell}(\G_0(L(N,D)),\psi^2)$.

\smallskip

Using somewhat extensive numerical experimentation, I have tried to see whether
the above bounds for $L(N,D)$ are optimal, and my conclusion is the following
conjecture, which I have considerably tested for quadratic characters, but
very little for non-quadratic ones, so I prefer stating it only for quadratic
characters, although it may well be true in general.

\begin{conjecture} Assume that $\psi$ is a quadratic character, that $k\ge2$,
  and that $D$ is a fundamental discriminant of sign $(-1)^k$.
  Set $\delta=4$ if $k$ is odd and $D\equiv1\pmod4$, and otherwise set
  $\delta=1$. We have  $g_{4,\delta D}\in M_{2k+4\ell}(\G_0(N),\psi^2)$,
  \emph{except} when $k$ is odd, $D\equiv12\pmod{16}$, and $N$ is odd
  (in which case we cannot say more than the theorem, i.e.,
  $g_{4,D}\in M_{2k+4\ell}(\G_0(2N),\psi^2)$).  
\end{conjecture}

\begin{corollary} Assume this conjecture. Then if $k$ is even, or $k$ is
  odd and $D\equiv0\pmod8$, we have $g_{4,D}\in M_{2k+4\ell}(\G_0(N),\psi^2)$.
\end{corollary}

\smallskip

\subsection{Application to $L(\chi_D,1-k)$ for $k$ Even}

\smallskip

The special case $k$ even, $\psi$ trivial, $N=1$, and $f=E_k$ (including
$k=2$) is well-known, at least for $\ell=0$:

\begin{corollary}\label{cor:keven} Let $k\ge2$ be an even integer.
  For any nonsquare integer $m$ define
  $$S_{k,\ell}(m)=\sum_{\substack{|s|<m^{1/2}\\s^2\equiv m\pmod 4}}P_{k,\ell}(m,s)\sigma_{k-1}\left(\dfrac{m-s^2}{4}\right)\;.$$
  For any positive fundamental discriminant $D$ let
  $$f_{k,\ell,D}(\tau)=-\delta_{\ell,0}\dfrac{B_k}{4k}L(\chi_D,1-k)+\sum_{n\ge1}\left(\sum_{d\mid n}\leg{D}{d}d^{k+2\ell-1}S_{k,\ell}((n/d)^2D)\right)q^n\;.$$

  We have $f_{k,\ell,D}\in M_{2k+4\ell}(\G)$, and it is a cusp form if $\ell>0$.
\end{corollary}

One of the important consequences is an explicit exact formula for
$L(\chi_D,1-k)$ for $D>0$ and $k$ even, due to Siegel. We begin by the
following:

\begin{theorem}\label{thm:sieg} Let $k$ be an even integer, let
  $r=\dim(M_k(\G))$, and define coefficients $c_i^{(k)}$ by
  $\Delta^{-r}E_{12r-k+2}=\sum_{i\ge-r}c_i^{(k)}q^i$ where by
  convention $E_0=1$. We have $c_0^{(k)}\ne0$, and for any
  $f=\sum_{n\ge0}a(n)q^n\in M_k(\G)$ we have
  $$\sum_{0\le n\le r}a(n)c_{-n}^{(k)}=0\;.$$
\end{theorem}

From this and the previous corollary, we deduce:

\begin{corollary}\label{cor:sieg}
  Let $k\ge2$ be an even integer, and let
  $r=\dim(M_{2k}(\G))=\lfloor k/6\rfloor+1$.
  With the above notation, for all positive fundamental discriminants $D$ we
  have
  $$L(\chi_D,1-k)=\dfrac{4k}{c_0^{(2k)}B_k}\sum_{1\le m\le r}S_{k,0}(m^2D)\sum_{1\le d\le r/m}\leg{D}{d}d^{k-1}c_{-dm}^{(2k)}\;.$$
\end{corollary}

\smallskip

\subsection{Application to $L(\chi_D,1-k)$ for $k$ Odd}

\smallskip

However a more general special case of the above theorem is important
to be able to treat the case $k$ odd.

First recall the following classical result on Eisenstein series:

\begin{proposition}
  For $s\in\C$ and $\psi$ a primitive Dirichlet character modulo $N$, define
  $$\sigma_s(\psi,n)=\sum_{d\mid n}\psi(d)d^s\;,$$
  and let $k\ge2$ be an integer such that $\psi(-1)=(-1)^k$.
  We have
  $$F_{k,\psi}(\tau):=\dfrac{L(\psi,1-k)}{2}+\sum_{n\ge1}\sigma_{k-1}(\psi,n)q^n\in M_k(\G_0(N),\psi)\;.$$
\end{proposition}

Note that there are slightly more general Eisenstein series involving
\emph{two} characters $\psi_1$ and $\psi_2$ and the more general function
$\sigma_s(\psi_1,\psi_2,n)=\sum_{d\mid n}\psi_1(d)\psi_2(n/d)d^s$, but we
will not need them.

\begin{corollary}
   Let $D$ be a positive fundamental discriminant, let
  $\psi$ be a primitive Dirichlet character modulo $N$, let $k\ge2$
  be an integer such that $\psi(-1)=(-1)^k$, and finally let $M$ be
  a positive integer. Set $\psi_1=\psi\chi_D\chi_{0,M}$, where as usual
  $\chi_{0,M}$ is the trivial character modulo $M$, and omitting most
  parameters from the notation, define
  \begin{align*}S_{k,\ell,M}(m)&=\sum_{\substack{|s|<m^{1/2}\\s^2\equiv m\pmod{M}}}P_{k,\ell}(m,s)\sigma_{k-1}\left(\psi,\dfrac{m-s^2}{M}\right)\text{\quad and}\\
    g_{M,D}(\tau)&=\delta_{\ell,0}\dfrac{L(\psi,1-k)L(\psi_1,1-k)}{4}+\sum_{n\ge1}\left(\sum_{d\mid n}\psi_1(d)d^{k-1}S_{k,\ell,M}((n/d)^2|D|)\right)q^n\;.\end{align*}
  We have $g_{4,D}\in M_{2k+4\ell}(\G_0(L(N,D)),\psi^2)$ with $L(N,D)$
  as in Definition \ref{def:lev}, we have
  $g_{1,D}\in M_{2k+4\ell}(\G_0(2N/\gcd(4,N)),\psi^2)$,
  and the forms $g_{M,D}$ are cusp forms when $\ell>0$.
\end{corollary}

\Proof Simply apply Theorem \ref{thm4} to
$F_{k,\psi}(M\tau)\in M_k(\G_0(MN),\psi\chi_{0,M})$. We have only included the
results that we need.\fp

\smallskip

One of the main applications of this corollary is the following: without
the auxiliary character $\psi$ and $N=1$, i.e., if we only use
the ordinary Eisenstein series of level $1$, in other words
Corollary \ref{cor:keven} with $\ell=0$, this allows to give very efficient
exact methods to compute $L(\chi_D,1-k)$ in terms of divisor sums
such as Corollary \ref{cor:sieg}, but only for $k$ even and $D>0$.

If we want to compute $L(\chi_D,1-k)$ for $k$ odd and $D<0$,
using the above corollary it is sufficient to use an auxiliary \emph{odd}
character $\psi$ such as $\psi=\chi_{-3}$ or $\chi_{-4}$, but again to
have a nonzero constant term we must choose $\ell=0$, hence
$P_{k,\ell}(m,s)=1$. We write explicitly a consequence of the special case
$\psi=\chi_{-4}$ (the case $\psi=\chi_{-3}$ would lead to forms of level
divisible by $6$, so much less interesting for practical use).

\begin{corollary} Let $D$ be a negative fundamental discriminant,
  let $t=|D|$ or $|D|/4$ be the corresponding positive
  squarefree integer, and let $k\ge3$ be an odd integer. For any nonsquare
  integer $m$, define
  $$S_k(m)=\sum_{|s|<m^{1/2}}\sigma_{k-1}(\chi_{-4},m-s^2)\;,\quad A_k=\left(1-\leg{D}{2}2^{k-1}\right)\dfrac{E_{k-1}}{8}\;,$$
  and set
  $$f_{k,D}(\tau)=A_kL(\chi_D,1-k)+\sum_{n\ge1}\left(\sum_{d\mid n}\leg{4D}{d}d^{k-1}S_k((n/d)^2t)\right)q^n\;,$$
  where the $E_k$ are the Euler numbers ($E_0=1$, $E_2=-1$, $E_4=5$, $E_6=-61$,...).

  Then for $D\ne-4$ we have $f_{k,D}\in M_{2k}(\G_0(2))$.
\end{corollary}

\Proof We apply the previous corollary with $\psi=\chi_{-4}$ hence $N=4$
and $L(\psi,1-k)=E_{k-1}/2$, with a fundamental discriminant $D_1$ equal
to the fundamental discriminant associated to $-4D$. We consider three cases:
if $D\equiv8\pmod{16}$ we have $D_1=-D$, hence $t=|D|/4$, and since by
the corollary we have $g_{4,D_1}\in M_{2k}(\G_0(2))$, the result follows
since it is clear that $g_{4,D_1}=f_{k,D}$, and
$\psi_1=\chi_{-4D_1}=\chi_{4D}=\chi_D$. If $D\equiv1\pmod4$ we have
$D_1=-4D$, hence $t=|D|$, and once again by the corollary we have
$g_{4,D_1}\in M_{2k}(\G_0(2))$ and it is clear that $g_{4,D_1}=f_{k,D}$.
However, in that case we have $\psi_1=\chi_{-4D_1}=\chi_{16D}=\chi_{4D}$,
so $L(\psi_1,1-k)=(1-\lgs{D}{2}2^{k-1})L(\chi_D,1-k)$.
Finally, if $D\equiv-4\pmod{16}$, we have $D_1=-D/4$ and $t=|D|/4$,
and by the corollary we have $g_{1,D_1}\in M_{2k}(\G_0(2))$ and we check
that $g_{1,D_1}=f_{k,D}$, and $\psi_1=\chi_{-4D_1}=\chi_D$. Finally,
we recall that $L(\chi_{-4},1-k)=E_{k-1}/2$.

Note that the restriction $D\ne-4$ is to ensure that $(n/d)^2t$ is never
a square. It would be easy to give a modification of the above corollary
for $D=-4$, but since the point is to give a formula for $L(\chi_D,1-k)$
this is unnecessary since we know that $L(\chi_{-4},1-k)=E_{k-1}/2$.\fp

\smallskip

To state the analog of Corollary \ref{cor:sieg} for $D<0$ and odd $k$,
we first introduce the following two modular forms in $M_4(\G_0(2))$:
$$F_4(\tau)=\dfrac{16E_4(2\tau)-E_4(\tau)}{15}\text{\quad and\quad}
\Delta_4(\tau)=\dfrac{E_4(\tau)-E_4(2\tau)}{240}\;.$$

Using an analog of Theorem \ref{thm:sieg} proved in Theorem 5.8 of
\cite{Coh4}, we deduce the following:

\begin{theorem}\label{thmsiegodd}
  Let $k\ge3$ be an odd integer, let $r=(k+3)/2$
  and define coefficients $c_i^{(k)}$ by
  $\Delta_4^{-r}F_4=\sum_{i\ge-r}c_i^{(k)}q^i$. We have $c_0^{(k)}\ne0$,
  and with the previous notation, for all negative fundamental discriminants
  $D\ne-4$ we have
  $$L(\chi_D,1-k)=-\dfrac{1}{c_0^kA_k}\sum_{1\le m\le r}S_k(m^2D)\sum_{1\le d\le r/m}\leg{D}{d}d^{k-1}c_{-dm}^{(k)}\;,$$
  where $A_k$ is as above.
\end{theorem}

\smallskip

We refer to \cite{Coh4} for applications of these types of results.

\section{Review of Hilbert Modular Forms}

We now come to the link with Hilbert modular forms (abbreviated HMF), and
review the basic definitions and notation.

Recall that if $K$ is a totally real number field of degree $d$, a Hilbert
modular form of (parallel) weight $k$ is a holomorphic function $F$ from the
product of $d$ copies of the upper half-plane $\H$ to $\C$ satisfying a
suitable modularity condition (holomorphy at the cusps is automatic for
$d\ge2$). To state the modularity condition, we need some notation. If
$a\in K$, we denote by $a^{(i)}$ for $1\le i\le d$ the $d$ embeddings of $a$
in $\R$. We define $\GL_2^+(\Z_K)$ as the group of
$\ga=\psmm{a}{b}{c}{d}\in M_2(\Z_K)$ such that $ad-bc$ is a totally positive
unit. For such $\ga$ we naturally denote by $\ga^{(i)}$ the element of
$\GL_2^+(\R)$ where the coefficients are sent to their $i$th embedding. We
make it act on $\H^d$ as follows: if $\bftau=(\tau_1,\dotsc,\tau_d)\in\H^d$
we define
$$\ga(\bftau)=\left(\ga^{(1)}(\tau_1),\dotsc,\ga^{(d)}(\tau_d)\right)\;,$$
where $\ga^{(i)}(\tau_i)=(a^{(i)}\tau_i+b^{(i)})/(c^{(i)}\tau_i+d^{(i)})$
is the natural action of $\GL_2^+(\R)$ on $\H$. It is also natural
to define the trace and the norm by
$$\Tr(\nu\bftau)=\sum_{1\le i\le d}\nu^{(i)}\tau_i\text{\quad and\quad}
\N(c\bftau+d)=\prod_{1\le i\le d}(c^{(i)}\tau_i+d^{(i)})\;.$$
The modularity condition can then simply be expressed as
$$F(\ga(\bftau))=\N(c\bftau+d)^kF(\bftau)$$
for all $\ga=\psmm{a}{b}{c}{d}\in\GL_2^+(\Z_K)$ (or suitable subgroups),
possibly with an additional multiplier system as in the one-variable case.

If $F$ is holomorphic and satisfies this, we say that it is a Hilbert modular
form of (parallel) weight $k$ (we will not use nonparallel weights). Such an
$F$ has a Fourier expansion
$$F(\bftau)=a_0+\sum_{\substack{\nu\in\gd^{-1}\\\nu\gg0}}a(\nu\gd)e^{2\pi i\Tr(\nu\bftau)}\;,$$
where the summation is over totally positive elements of the codifferent of
$K$. Note that the Fourier coefficients depend only on the \emph{ideal}
$\nu\gd$ (in other words are the same for $\nu$ and $\nu\eps$ for all units
$\eps$), but if we had chosen the slightly weaker notion of modularity
under $\SL_2(\Z_K)$, this would have been the case only when totally positive
units of $K$ are all squares.

Our goal is to give explicit formulas showing how to obtain one-variable
(``elliptic'') modular forms from Hilbert modular forms and conversely.
As already mentioned, for consistency, we will use lower-case letters such as
$f$, $g$, $h$ for one-variable functions, and upper-case letters such as $F$
for functions of several variables.

\section{Obtaining Elliptic Modular Forms from Hilbert Modular Forms}

We now study the (easy) methods to obtain one-variable modular forms from HMFs.
 
The simplest and most evident way is to restrict to the diagonal: if $F$ is a
HMF of parallel weight $k$ on a totally real number field of degree $d$, the
function $f(\tau)=F(\tau,\tau,\dots,\tau)$ is trivially
a one-variable modular form of weight $kd$. But one can obtain many more
such functions by using suitable \emph{partial derivatives} of $F$: the
partial derivatives themselves are modular with parasitic terms,
but certain linear combinations make these terms disappear. We have
the following result:

\begin{theorem}\label{thm1} Let $F$ be a HMF of weight $k$ on a totally real
field of degree $d$, let $U=(u_1,u_2,\dotsc,u_d)\in\C^d$ be a vector of $d$
complex numbers such that $\sum_{1\le i\le d}u_i=0$, and let $\ell$ be a
nonnegative integer. To simplify notation, denote by $\partial_{\tau_i}$ the
partial derivative with respect to $\tau_i$. The function
$$D_{U,\ell}(F)(\tau)=\sum_{\substack{r_1,r_2,\dotsc,r_d\ge0\\r_1+r_2+\cdots r_d=\ell}}\prod_{1\le i\le d}\dfrac{u_i^{r_i}(k-1)!}{(k-1+r_i)!}\dfrac{\partial_{\tau_i}^{r_i}}{r_i!}(F)(\tau,\tau,\dotsc,\tau)$$
is a modular form of weight $dk+2\ell$, and a cusp form if $\ell>0$.
\end{theorem}

\Proof Let $\bfz=(z_1,z_2,\dotsc,z_d)$ be a vector of $d$ formal variables,
and generalizing the notation of Corollary \ref{cor:FG}, set
$$B_F(\bftau,\bfz)=\sum_{r_1,r_2,\dotsc,r_d\ge0}\prod_{1\le i\le d}\dfrac{(2\pi iz_i^2)^{r_i}}{r_i!(k-1+r_i)!}\partial_{\tau_i}^{r_i}(F)(\bftau)\;.$$
By the same calculation as in the proof of Corollary \ref{cor:FG},
 if $\ga=\psmm{a}{b}{c}{d}\in\GL_2^+(\Z_K)$
we have, with obvious termwise action on $\bfz$:
$$B_F|_k\ga(\bftau,\bfz)=e^{2\pi iS({\bftau},\bfz)}B_{F|_k\ga}(\bftau,\bfz)\;,\text{\ with\ }S(\bftau,\bfz)=\sum_{1\le i\le d}z_i^2c^{(i)}/(c^{(i)}\tau_i+d^{(i)})\;.$$
Thus, in particular if we restrict to $\psmm{a}{b}{c}{d}\in\SL_2(\Z)$
and $\bftau$ to the diagonal $\bftau=(\tau,\tau,\dotsc,\tau)$, we have
$S(\bftau,\bfz)=(c/(c\tau+d))\sum_{1\le i\le d}z_i^2$. If as in the
theorem we choose $z_i=u_i^{1/2}z$ with $\sum_{1\le i\le d}u_i=0$, we have
$S(\bftau,\bfz)=0$, proving that $B_F$ restricted to the diagonal behaves
modularly if $F$ does. Since this is now a formal power series in the single
variable $z$, this means that each coefficient is modular, proving the
theorem since the coefficient of $z^{2\ell}$ is the sum of all terms such
that $r_1+r_2+\cdots r_d=\ell$.

We will see below that the converse of this theorem is essentially true.

\smallskip

To compare with the results of the next section, we specialize to the
case $u_1=-1$, $u_2=1$ and $u_i=0$ for $i\ge3$:

\begin{corollary}\label{cor:Dell}
  Let $F$ be a HMF of weight $k$ on a totally real number
  field of degree $d$, and let $\ell$ be a nonnegative integer. The function
$$D_{\ell}(F)(\tau)=\sum_{0\le r\le\ell}(-1)^r\dfrac{(k-1)!^2}{(k-1+r)!(k-1+\ell-r)!}\dfrac{\partial_{\tau_1}^r\partial_{\tau_2}^{\ell-r}}{r!(\ell-r)!}(F)(\tau,\tau,\dotsc,\tau)$$
is a modular form of weight $dk+2\ell$, and a cusp form if $\ell>0$.
\end{corollary}

\section{Obtaining Elliptic Modular Forms from Hecke--Eisenstein Series}

\smallskip

\subsection{The General Formula}

\smallskip

As in the one-variable case, Eisenstein series (here called
Hecke--Eisenstein series) give a supply of HMF with explicit Fourier
expansion, and applying the results of the previous section we can thus
obtain interesting one-variable forms. The formula used, due to Siegel, is
the following:

\begin{theorem}\label{thmeisk} Let $K$ be a totally real number field of
  degree $d$, and let $k\ge2$ be an even integer. For any ideal $I$ and complex
  number $s$, define the ideal sum of divisors function by
  $\sigma_{K,s}(I)=\sum_{J\mid I}\N(J)^s$.
  The function
  $$\dfrac{\z_K(1-k)}{2^d}+\sum_{\substack{\nu\in\delta^{-1}\\\nu\gg0}}\sigma_{K,k-1}(\nu\delta)e^{2\pi i\Tr(\nu\bftau)}$$
    is a Hilbert modular form of weight $k$ on the full group $\GL_2^+(\Z_K)$.
\end{theorem}

Clearly
$$\partial^{r_i}_{\tau_i}(e^{2\pi i\Tr(\nu\bftau)})=(2\pi i\nu^{(i)})^{r_i}e^{2\pi i\Tr(\nu\bftau)}\;,$$
so we deduce from Theorem \ref{thm1} the following:

\begin{theorem}\label{thmeiskspec} Let $U=(u_1,\dotsc,u_d)\in\C^d$ be a vector
of $d$ complex numbers such that $\sum_{1\le i\le d}u_i=0$, and let $\ell$ be a
nonnegative integer. For $\nu\in K$ define
$$c_{U,\ell}(\nu)=\sum_{\substack{r_1,\dotsc,r_d\ge0\\r_1+\cdots+r_d=\ell}}\prod_{1\le i\le d}\dfrac{(k-1)!}{r_i!(k-1+r_i)!}(u_i\nu^{(i)})^{r_i}\;.$$
The function
 $$\delta_{\ell,0}\dfrac{\z_K(1-k)}{2^d}+\sum_{n\ge1}\left(\sum_{\substack{\nu\in\delta^{-1},\ \nu\gg0\\\Tr(\nu)=n}}c_{U,\ell}(\nu)\sigma_{K,k-1}(\nu\delta)\right)q^n$$
belongs to $M_{dk+2l}(\G)$, and is a cusp form if $\ell>0$.
\end{theorem}

From this for $\ell=0$ and Theorem \ref{thm:sieg} we deduce a finite
explicit formula for $\z_K(1-k)$:

\begin{corollary}\label{cor:zk} Let $k\ge2$ be an even integer, and let
  $r=\dim(M_{dk}(\G))$. We have
  $$\z_K(1-k)=-\dfrac{2^d}{c_0^{(dk)}}\sum_{1\le n\le r}c_{-n}^{(dk)}\sum_{\substack{\nu\in\gd^{-1},\ \nu\gg0\\\Tr(\nu)=n}}\sigma_{K,k-1}(\nu\gd)\;,$$
  where the $c_{-n}^{(dk)}$ are given in Theorem \ref{thm:sieg}.
\end{corollary}

\smallskip

\subsection{Practical Application}

\smallskip

To use the above formulas in practice, we need to be able to describe
explicitly the set
$$E_n=\{\nu\in\gd^{-1},\ \nu\gg0,\ \Tr(\nu)=n\}\;.$$
This can be done as follows. We let $(1,u_2,\dotsc,u_d)$ be an integral
basis, i.e., a $\Z$-basis of $\Z_K$ chosen so that the first element $u_1$ is
equal to $1$. Choose as basis of the codifferent the dual basis
$(v_1,v_2,\dotsc,v_d)$ obtained by inverting the matrix
$(\Tr(u_iu_j))_{1\le i,j\le d}$. Finally, denote by
$(\sigma_1,\sigma_2,\dotsc,\sigma_d)$ the $d$ embeddings of $K$ into $\R$.

\begin{proposition}\label{prop:simplex} Denote by $\Delta$ the simplex of
  $\R^{d-1}$ with vertices
  $$P_i=(\sigma_i(u_2),\sigma_i(u_3),\dotsc,\sigma_i(u_d))_{1\le i\le d}\;,$$
  and let $\Delta_n$ be the set of points in the interior of $\Delta$
  whose coordinates are integral multiples of $1/n$. We have
  $$E_n=\{n(v_1+x_2v_2+x_3v_3+\cdots+x_dv_d),\ (x_2,x_3,\dotsc,x_d)\in \Delta_n\}\;.$$
\end{proposition}

\Proof By definition, we have $\Tr(u_iv_j)=\delta_{i,j}$.
Any element $\nu\in\gd^{-1}$ is of the form $\nu=\sum_{1\le i\le d}y_iv_i$
with $y_i\in\Z$, hence since $\Tr(v_1)=1$ and $\Tr(v_i)=0$ for $i>1$, we have
$\Tr(\nu)=y_1$, hence if $\nu\in E_n$ we have $y_1=n$, so setting $x_i=y_i/n$
for $i\ge2$ we have $\nu=n(v_1+\sum_{2\le i\le d}x_iv_i)$, with the $x_i$
integral multiples of $1/n$. It is easy to check that the
$(\sigma_i(\nu)/n)_{1\le i\le d}$ are the barycentric coordinates of
$(x_2,\dotsc,x_d)$ with respect to the vertices of $\Delta_n$, so that
$\nu\gg0$ is equivalent to $(x_2,\dotsc,x_d)$ being in the interior of
$\Delta$.\fp

Computationally, the formula of Corollary \ref{cor:zk} can be used in two
ways. First, as a direct way to compute $\z_K(1-k)$, using the above
proposition to range in the finite set of totally positive $\nu\in\gd^{-1}$
of given trace, but of course this is slow for large degrees. But
more importantly, since it gives a trivial multiple of the \emph{denominator}
of the rational number $\z_K(1-k)$, we can use one of several \emph{numerical}
methods to approximate $\z_K(1-k)$, and use error bounds plus the bound
on the denominator to obtain the exact value, and this is almost always
considerably faster.

\smallskip

\subsection{The Case of Real Quadratic Fields}

\smallskip

In the case where $K$ is a real quadratic field we can even be more
explicit since the elements of $\Z_K$ and of $\delta$ can trivially be
written down. We first specialize Theorem \ref{thmeisk}:

\begin{theorem} Let $K$ be a real quadratic field of discriminant $D$
  and $k\ge2$ be an even integer. Define
  $$c_k(s,n)=\sum_{\substack{d\mid\gcd(s,n)\\(s/d)\equiv(n/d)D\pmod2}}\leg{D}{d}d^{k-1}\sigma_{k-1}\left(\dfrac{(n/d)^2D-(s/d)^2}{4}\right)\;.$$
  The function
  $$\dfrac{\z_K(1-k)}{4}+\sum_{n\ge1}\sum_{\substack{|s|<n\sqrt{D}\\s\equiv nD\pmod2}}c_k(s,n)e^{2\pi i\Tr(((n+s/\sqrt{D})/2)\bftau)}$$
  is a Hilbert modular form of weight $k$ on $\GL_2^+(\Z_K)$.
\end{theorem}

\Proof Let $K=\Q(\sqrt{D})$ with $D>0$ a fundamental discriminant. The
elements $\nu$ of the codifferent are the elements of $K$ of the form
$\nu=(n+s/\sqrt{D})/2$ with $s\equiv nD\pmod2$, and the condition
$\nu\gg0$ means that $|s|<n\sqrt{D}$ (the condition $\Tr(\nu)=n$ is implicit
in the notation). The ideal $\nu\delta$ is the ideal $((s+n\sqrt{D})/2)\Z_K$,
so the main work is the computation of $\sigma_{K,k-1}(\nu\delta)$.
By definition we have
$$\sigma_{K,k-1}(I)=\sum_{d\mid\N(I)}d^{k-1}\sum_{\substack{J\mid I\\\N(J)=d}}1\;,$$
so we must compute this inner sum for $I=((s+n\sqrt{D})/2)\Z_K$, hence
$\N(I)=(n^2D-s^2)/4$. We have the following lemma:

\begin{lemma}\label{lem:sum1} With the above notation, we have
  $$\sum_{\substack{J\mid I\\\N(J)=d}}1=\sum_{e\mid\gcd(s,n,d,(n^2D-s^2)/(4d))}\leg{D}{e}\;.$$
\end{lemma}

\Proof The proof of this lemma is of a combinatorial nature: we first note
that both sides are multiplicative functions of $d$, and then for $d=p^{\al}$
for a prime $p$ we distinguish between the cases $p$ inert, split, and
ramified. Since this proof (essentially due to Siegel) is classical and tedious
we omit the details.\fp

Once this lemma is proved, it is a simple matter of bookkeeping to obtain
the theorem, and the details are left to the reader.\fp

\smallskip

It is now immediate to specialize Theorem \ref{thmeiskspec} to the
quadratic case, and in fact up to a multiplicative constant without loss of
generality we can choose $u_1=-1$ and $u_2=1$. After a small computation
we obtain the following result whose detailed proof is left to the reader:

\begin{corollary}\label{cor:qk} Keep the same notation as the theorem. Define
  $$Q_{k,\ell}(n,s,D)=\sum_{0\le r\le \ell}(-1)^r\dfrac{(k-1)!^2}{(k-1+r)!(k-1+\ell-r)!}\dfrac{(n+s/\sqrt{D})^r(n-s/\sqrt{D})^{\ell-r}}{r!(\ell-r)!}$$
  and
  $$S_{k,\ell}(m,D)=\sum_{\substack{|s|<m\sqrt{D}\\s\equiv mD\pmod2}}Q_{k,\ell}(m,s,D)\sigma_{k-1}\left(\dfrac{m^2D-s^2}{4}\right)\;.$$
  The function
$$\delta_{\ell,0}\dfrac{\z_K(1-k)}{4}+\sum_{n\ge1}\sum_{d\mid n}\leg{D}{d}d^{k+\ell-1}S_{k,\ell}(n/d,D)q^n$$
  belongs to $M_{2k+2\ell}(\G)$ and is a cusp form if $\ell>0$.
\end{corollary}

Since for $\ell$ odd $Q_{k,\ell}(n,s,D)$ is an odd polynomial in $s$ and since
the summation on $s$ is symmetrical, it follows that for $\ell$ odd the form is
identically zero, so we can replace $\ell$ by $2\ell$ in which case the
formula of this corollary is, up to a constant multiplicative factor,
identical to that of Corollary \ref{cor:keven}. More precisely we leave to
the reader the proof of the following result, already stated in Proposition
\ref{prop:pkl}:

\begin{proposition}\label{prop:pq} We have
$$P_{k,\ell}(n^2D,s)=(-1)^\ell\dfrac{(2\ell)!}{2^{2\ell}\ell!}\dfrac{(k+2\ell-1)!(k+\ell-1)!}{(k-1)!^2}D^\ell Q_{k,2\ell}(n,s,D)\;.$$
\end{proposition}

\section{Obtaining Hilbert Modular Forms from Elliptic Modular Forms}

Probably the prototypical ``lifting theorem'' from one-variable modular
forms to Hilbert modular forms is due to Doi--Naganuma \cite{Doi-Nag}.
I will give here a slightly different version which is both more general,
in that it allows arbitrary modular forms, and weaker, since the level is
in general not optimal, and no Hecke equivariant statement is given. First,
we prove the following partial converse to Theorem \ref{thm1}:

\begin{theorem}\label{thm1conv}
  Let $F$ be a holomorphic function on $\H^2$, and assume that
  all the functions $D_\ell(F)$ of Corollary \ref{cor:Dell} belong to
  $M_{2k+2\ell}(\G_0(N),\psi^2)$ for some $N$ and $\psi$.
  \begin{enumerate}\item The function $F$ is modular of parallel weight
    $(k,k)$ on $\G_0(N)$ considered as a subgroup of $\GL_2^+(\Z_K)$ with
    character $\psi\circ\N$.
\item If, in addition, $F$ has a Fourier expansion of the standard form
$$F(\bftau)=a_0+\sum_{\substack{\nu\in\gd^{-1}\\\nu\gg0}}a(\nu\gd)e^{2\pi i\Tr(\nu\bftau)}\;,$$
  then it is modular (with the same character) on the subgroup
  $\widetilde{G}$ generated by $\G_0(N)$ and
  the matrices $\psmm{\eps_1}{b}{0}{\eps_2}\in\GL_2^+(\Z_K)$ with $\eps_1$ and
  $\eps_2$ units of $\Z_K$, which is a subgroup of finite index of
  $\GL_2^+(\Z_K)$.
\item If $D\equiv1\pmod4$ and $F$ has the above Fourier expansion, it is
  modular on the subgroup $\G_0(N\Z_K)$ of matrices
  $\psmm{a}{b}{c}{d}\in\GL_2^+(\Z_K)$ such that $c\in N\Z_K$.
\end{enumerate}
\end{theorem}

\Proof (1). We must first prove that if
$\ga=\psmm{a}{b}{c}{d}\in\G_0(N)$ we have $F|k\ga=\psi(d)^2F$, in other words
that if we set
$$h(\tau_1,\tau_2)=F\left(\dfrac{a\tau_1+b}{c\tau_1+d},\dfrac{a\tau_2+b}{c\tau_2+d}\right)-\psi(d)^2(c\tau_1+d)^k(c\tau_2+d)^kF(\tau_1,\tau_2)\;,$$
we have $h=0$. Since $h$ is a holomorphic function in two variables, it
suffices to prove that all the partial derivatives $\partial_{\tau_1}^{r_1}\partial_{\tau_2}^{r_2}(h)(\tau,\tau)$ vanish. We will prove this
by induction on $\ell=r_1+r_2$. First, note that, using the function $B_F$
defined above, the modularity of all the $D_{\ell}(F)$ is equivalent to
the single identity
$$B_F\left(\dfrac{a\tau+b}{c\tau+d},\dfrac{a\tau+b}{c\tau+d},\dfrac{z}{c\tau+d},\dfrac{iz}{c\tau+d}\right)=\psi(d)^2(c\tau+d)^{2k}B_F(\tau,\tau,z,iz)$$
as a formal power series in $z$ (this is in fact how the $D_{\ell}(F)$
were constructed).

By Corollary \ref{cor:FG} (more precisely its generalization to two variables)
we have $B_F|_k\ga(\bftau,\bfz)=e^{S(\bftau,\bfz)}B_{F|_k\ga}(\bftau,\bfz)$.
Combining this with the above identity implies that
$$\sum_{r_1,r_2\ge0}\dfrac{z^{2r_1}(iz)^{2r_2}}{(k-1+r_1)!(k-1+r_2)!}\dfrac{\partial_{\tau_1}^{r_1}\partial_{\tau_2}^{r_2}(h)(\tau,\tau)}{r_1!r_2!}=0\;,$$
so isolating the term in $z^{2\ell}$ we deduce that for all $\ell\ge0$ and
  $$\sum_{0\le r\le\ell}(-1)^r\dfrac{\partial_{\tau_1}^r\partial_{\tau_2}^{\ell-r}(h)(\tau,\tau)}{(k-1+r)!(k-1+\ell-r)!r!(\ell-r)!}=0\;.$$
  Now assume by induction on $\ell$ that for all $r$ with $0\le r\le \ell-1$
  we have $\partial_{\tau_1}^r\partial_{\tau_2}^{\ell-1-r}(h)(\tau,\tau)=0$,
  which is true for $\ell=1$ by the above identity. Computing the derivative
  of this with respect to $\tau$ we deduce that
  $(\partial_{\tau_1}^{r+1}\partial_{\tau_2}^{\ell-1-r}+\partial_{\tau_1}^r\partial_{\tau_2}^{\ell-r})(h)(\tau,\tau)=0$, so by induction we deduce that
  $$\partial_{\tau_1}^r\partial_{\tau_2}^{\ell-r}(h)(\tau,\tau)=(-1)^{\ell-r}\partial_{\tau_1}^\ell(h)(\tau,\tau)=(-1)^r\partial_{\tau_2}^\ell(h)(\tau,\tau)\;.$$
  Replacing in the above identity we obtain
  $$\partial_{\tau_2}^\ell(h)(\tau,\tau)\sum_{0\le r\le \ell}\dfrac{1}{(k-1+r)!(k-1+\ell-r)!r!(\ell-r)!}=0\;,$$
  hence $\partial_{\tau_2}^\ell(h)(\tau,\tau)=0$, hence all the
  $\partial_{\tau_1}^r\partial_{\tau_2}^{\ell-r}(h)(\tau,\tau)$ vanish, proving
  our induction hence (1).

  \medskip

  (2). The modularity for $\widetilde{G}$ is immediate from the assumption
  by multiplicativity. The fact that this is of finite index is not difficult
  but we omit the proof.

  \medskip

  (3). This is proved in \cite{Lio-Ver}. It is possible that the result is
  also true for $D\equiv0\pmod4$.\fp

  Note that by using the same methods, if desired one can also prove a
  converse to the more general Theorem \ref{thm1}.

  \medskip
  
\begin{theorem}\label{thm3} Let $f=\sum_{n\ge0}a(n)q^n\in M_k(\G_0(N),\psi)$
with $k$ even, let $K=\Q(\sqrt{D})$ be a real quadratic field of discriminant
$D$, and let $\psi_1=\psi\chi_D$. Let
$$A_k(\nu\gd)=\sum_{d\mid \nu\gd,\ d\in\Z_{>0}}d^{k-1}\psi_1(d)a(\N(\nu\gd/d))\;.$$
The two-variable function on $\H^2$:
$$F(\bftau)=a(0)\dfrac{L(\psi_1,1-k)}{2}+\sum_{\substack{\nu\in\gd^{-1}\\\nu\gg0}}e^{2\pi i\Tr(\nu\bftau)}A_k(\nu\gd)\;,$$
where $\bftau=(\tau_1,\tau_2)$, is a Hilbert modular form of parallel weight
$k$ on some subgroup $\widetilde{G}$ of $\GL_2^+(\Z_K)$ of level dividing $2N$
(more precisely dividing $L(N,D)$ defined in Definition \ref{def:lev})
with character $\psi\circ\N$, in other words for all
$\ga=\psmm{a}{b}{c}{d}\in \widetilde{G}$
we have $F|_k\ga=\psi(\N(d))F$.

In addition, if $D\equiv1\pmod4$, $F$ is modular on the subgroup
 $\G_0(N\Z_K)$ of matrices $\psmm{a}{b}{c}{d}\in\GL_2^+(\Z_K)$ such that
 $c\in N\Z_K$.
\end{theorem}

Note that $d\mid\nu\gd$ means that $d\Z_K\mid\nu\gd$, in other words
$\nu\gd\subset d\Z_K$.

\smallskip

\Proof Recall from Lemma \ref{lem:sum1} that the totally positive elements
of $\gd^{-1}$ of trace $n$ are of the form $\nu=(n+s/\sqrt{D})/2$ with
$|s|<n\sqrt{D}$ and $s\equiv nD\pmod2$. Thus, for such a $\nu$ we have
$$A_k(\nu\delta)=\sum_{\substack{d\mid\gcd(n,s)\\s/d\equiv (n/d)D\pmod2}}d^{k-1}\psi_1(d)a\left(\dfrac{(n/d)^2D-(s/d)^2}{4}\right)\;.$$
Thus
$$F(\tau,\tau)=a(0)\dfrac{L(\psi_1,1-k)}{2}+\sum_{n\ge1}\left(\sum_{d\mid n}d^{k-1}\psi_1(d)\sum_{\substack{|s|<(n/d)\sqrt{D}\\s\equiv(n/d)D\pmod2}}a\left(\dfrac{(n/d)^2D-s^2}{4}\right)\right)q^n\;,$$
which is exactly the function $g_{4,D}$ of Theorem \ref{thm4} for $\ell=0$, so
$F(\tau,\tau)\in M_{2k}(\G_0(L(N,D)),\psi^2)$.

\smallskip

More generally, for $\ell>0$, since
$\partial^{r_i}_{\tau_i}(e^{2\pi i\Tr(\nu\bftau)})=(2\pi i\nu^{(i)})^{r_i}e^{2\pi i\Tr(\nu\bftau)}$, we compute that with the notation of Corollary
\ref{cor:Dell} we have
$$\dfrac{D_{\ell}(F)(\tau)}{(\pi i)^{\ell}}=\sum_{n\ge1}\left(\sum_{d\mid n}d^{k+\ell-1}\psi_1(d)\sum_{\substack{|s|<(n/d)\sqrt{D}\\s\equiv(n/d)D\pmod2}}Q_{k,\ell}(n/d,s,D)a\left(\dfrac{(n/d)^2D-s^2}{4}\right)\right)q^n\;,$$
with $Q_{k,\ell}$ as in Corollary \ref{cor:qk}. As before, since the sum on $s$
is symmetrical we have $D_{\ell}(F)=0$ for $\ell$ odd, and using the relation
given above between $Q_{k,2\ell}$ and $P_{k,\ell}$ we deduce that
$$D_{2\ell}(F)(\tau)=C\sum_{n\ge1}\left(\sum_{d\mid n}d^{k+2\ell-1}\psi_1(d)\sum_{\substack{|s|<(n/d)\sqrt{D}\\s\equiv(n/d)D\pmod2}}P_{k,\ell}((n/d)^2D,s)a\left(\dfrac{(n/d)^2D-s^2}{4}\right)\right)q^n\;,$$
for some irrelevant constant $C$, which is again a constant multiple
of $g_{4,D}$ of Theorem \ref{thm4} for $\ell$. We deduce that for all
$\ell$ we have $D_{\ell}(F)\in M_{2k}(\G_0(L(N,D)),\psi^2)$,
hence our theorem follows from Theorem \ref{thm1conv}.\fp

\medskip

{\bf Remark.} All our proofs based on handling of explicit Fourier
coefficients are in some sense ``naive''. The best way to prove the kind
of classical results given in the present paper is to use suitable
\emph{kernel functions}, which give the same formulas, but with better
control on the level. For instance this is how the quoted results of
\cite{Lio-Ver} are proved.

\end{document}